\documentclass[11pt]{article}
\usepackage[T1]{fontenc}
\usepackage[utf8]{inputenc}
\usepackage{lmodern,microtype}
\usepackage[a4paper,margin=28mm]{geometry}
\usepackage{amsmath,amssymb,amsthm,mathtools}
\usepackage{booktabs,array,enumitem,needspace}
\usepackage[colorlinks=true,linkcolor=blue!60!black,citecolor=blue!60!black,urlcolor=blue!60!black]{hyperref}
\usepackage{xcolor}
\hypersetup{pdftitle={Applications of reduced and coreduced modules III: Derived torsion, completion, and coefficientwise vanishing},pdfauthor={David Ssevviiri}}
\usepackage{fancyhdr}
\newtheorem{theorem}{Theorem}[section]
\newtheorem{proposition}[theorem]{Proposition}
\newtheorem{corollary}[theorem]{Corollary}
\newtheorem{lemma}[theorem]{Lemma}
\theoremstyle{definition}
\newtheorem{definition}[theorem]{Definition}
\newtheorem{example}[theorem]{Example}
\newtheorem{question}[theorem]{Question}
\theoremstyle{remark}\newtheorem{remark}[theorem]{Remark}
\DeclareMathOperator{\Hom}{Hom}\DeclareMathOperator{\Ext}{Ext}
\DeclareMathOperator{\Tor}{Tor}\DeclareMathOperator{\ann}{ann}
\DeclareMathOperator{\coker}{coker}\DeclareMathOperator{\im}{im}
\DeclareMathOperator{\pd}{pd}\DeclareMathOperator{\fd}{fd}
\DeclareMathOperator{\cdim}{cd}\DeclareMathOperator{\hdim}{hd}

\newcommand{\Mod}{\text{-}\mathrm{Mod}}
\newcommand{\modf}{\text{-}\mathrm{mod}}
\newcommand{\NN}{\mathbb N}
\newcommand{\LL}{\mathrm L}
\newcommand{\RR}{\mathrm R}
\newcommand{\colim}{\varinjlim}
\newcommand{\ilim}{\varprojlim}
\title{Applications of reduced and coreduced modules III: \\Derived torsion and completion over polynomial extensions}
\author{David Ssevviiri\\[2pt]\small Department of Mathematics, Makerere University\\\small P.O. Box 7062, Kampala, Uganda\\\small\texttt{david.ssevviiri@mak.ac.ug}}
\date{}
\begin{document}
\maketitle
\begin{abstract}
We study when degree-zero identities for reduced and coreduced modules extend to ambient derived torsion and completion. For a non-idempotent ideal $I$ in an explicit polynomial ring, we construct a module $M$ that is simultaneously $I$-reduced and $I$-coreduced and satisfies $\Gamma_I(M)=\Lambda_I(M)=0$ but $H_I^2(M)\ne0$. The module is a union of finite reduced and coreduced submodules whose positive derived torsion vanishes. The supporting ideal has unrestricted torsion dimension two, whereas the dimension on each bounded-nilpotence quotient category is one. The calculations use compatible cone resolutions for $R=S[t]$ and $I=JS[t]+(t)$, where $J^2=J$. These give paired Ext--Tor formulas for coefficients annihilated by a fixed power of $t$. For arbitrary coefficients, an explicit mapping-fiber model realizes the standard derived-composition comparison sequence and identifies the obstruction to extending the bounded formula. Coordinate calculations for countably idempotent-generated ideals supply the computational foundation, and a valuation-domain example gives nonzero first derived completion on a cyclic reduced and coreduced coefficient.
\end{abstract}

\section{Introduction}
If an $R$-module $M$ is $I$-reduced or $I$-coreduced in the ideal-power sense, the corresponding degree-zero functor simplifies to
\[
 \Gamma_I(M)\cong\Hom_R(R/I,M),\qquad
 \Lambda_I(M)\cong(R/I)\otimes_RM,
\]
respectively; see \cite{SsevI,SsevII} and Lemma~\ref{lem:degreezero}. These identities concern the coefficient module itself, whereas ambient derived functors use resolutions in the entire module category. The question guiding this paper is:
\begin{quote}
When do the degree-zero simplifications associated with reduced and coreduced modules extend to computations of ambient derived torsion and completion?
\end{quote}
We answer this question through explicit comparisons and a counterexample over a polynomial extension. The counterexample is the principal module-theoretic application.

\Needspace{22\baselineskip}
\begin{theorem}[Main application]\label{thm:main}
Let $k$ be a field, and put
\[
 S=\prod_{n\ge1}k,\qquad J=\bigoplus_{n\ge1}k,\qquad
 R=S[t],\qquad I=JR+(t).
\]
Give $V_n=k[t]/(t^n)$ the $S$-action through the $n$th coordinate, and set $M=\bigoplus_{n\ge1}V_n$. Then $I^2\ne I$, the module $M$ is both $I$-reduced and $I$-coreduced, and
\[
 \Gamma_I(M)=\Lambda_I(M)=0,\qquad H_I^2(M)\ne0.
\]
For every $N\ge1$, the submodule $M_N=\bigoplus_{1\le n\le N}V_n$ is finitely generated, is $I$-reduced and $I$-coreduced, and satisfies $H_I^q(M_N)=0$ for $q>0$. Nevertheless $M=\colim_NM_N$ has the displayed nonvanishing. Moreover,
\[
 \operatorname{cd}_R(I)=2,\qquad
 \sup\{q\ge0:H_I^q(X)\ne0\text{ for some }X\text{ with }t^NX=0\}=1
 \quad(N\ge1).
\]
\end{theorem}
\begin{proof}[Proof of the main application]
The coefficient calculation and the assertions about $M_N$ are proved in Theorem~\ref{ex:localnil}; the dimension statements are Corollary~\ref{cor:dimensionjump}.
\end{proof}
The hypothesis on the coefficient is substantive here: the supporting ideal is not idempotent. The example shows that stabilization in degree zero, even for a module satisfying both conditions simultaneously, does not imply ambient torsion-acyclicity. It also distinguishes a fixed bound on the nilpotence of $t$ from local nilpotence on individual elements.

\subsection*{Comparison results underlying the example}
For a commutative ring $S$ and an idempotent ideal $J$, put $T=S/J$, $R=S[t]$, and $I=JR+(t)$. Theorem~\ref{thm:bounded} computes, for $t^NM=0$,
\[
 H_I^q(M)\cong\Ext_S^q(T,M),\qquad
 \LL_q\Lambda_I(M)\cong\Tor_q^S(T,M).
\]
It includes the transition maps at every sufficiently large ideal-power stage, and identifies both calculations with the relative derived functors in $S[t]/(t^N)\Mod$. The torsion half for $N=1$ is shared with the companion manuscript \cite[Theorem 5.1]{MS}; the paired calculation here treats every fixed $N$.

For arbitrary coefficients and $q\ge1$, write $E^q=\Ext_S^q(T,M)$, with its $t$-action. Theorem~\ref{thm:fiber} gives the explicit comparison sequence
\[
 0\longrightarrow E^{q-1}[t^{-1}]/\im(E^{q-1})
 \longrightarrow H_I^q(M)
 \longrightarrow\ker\bigl(E^q\to E^q[t^{-1}]\bigr)
 \longrightarrow0.
\]
This is a polynomial realization of the standard derived-composition construction described below. Its terms make the failure in Theorem~\ref{thm:main} explicit: local nilpotence of $t$ on $M$ need not pass to $\Ext_S^1(T,M)$. The valuation-domain calculation in Example~\ref{ex:valuation} supplies the complementary phenomenon of nonzero first projective-derived completion on a cyclic reduced and coreduced coefficient.

\subsection*{Classical inputs, existing comparisons, and overlap}
The roles of the relevant earlier results are as follows.
\begin{enumerate}[label=(\roman*)]
\item The ideal-power degree-zero identities and the global equivalence between idempotence and reducedness or coreducedness of every module are inputs from \cite{SsevI,SsevII}. Classical projectivity of countably idempotent-generated ideals \cite[Tag 05WH]{Stacks} gives the coordinate calculations in Section~\ref{sec:coordinates}. Their finite-coefficient and sum criteria are elementary consequences of Hom--Ext arguments; the broader theory of preservation of sums by defect and Ext functors is developed in \cite{BZ}.
\item Qu\'y and Rohrer \cite[Proposition 3.4]{QR} give a comparison sequence for adjoining a principal ideal under hypotheses that the relevant torsion functors preserve injectives. Their sequence has the same kernel--cokernel form as the sequence displayed above. In Section~\ref{sec:obstruction} we verify directly the weaker acyclicity needed for the composite $\Gamma_{(t)}\circ\Hom_R(T[t],-)$, for arbitrary idempotent $J$, and supply its compatible cone model. The derived-composition and spectral-sequence tools are standard \cite[Tags 015M, 015N]{Stacks}.
\item Positive cohomology for infinitely generated ideals over non-Noetherian absolutely flat rings is already established in \cite[Proposition 3.7]{QR}. The cyclic product-of-fields calculation used in Section~\ref{sec:coordinates} also occurs in \cite[Example 5.4]{MS}. Here it is a starting calculation for the polynomial coefficient in Theorem~\ref{thm:main}, whose degree-two value and reduced/coreduced properties are checked explicitly.
\item Weak proregularity provides comparisons with \v{C}ech and sequential constructions in the finitely generated setting \cite{PSY,PSYerr,VY,Yek}. The present calculations use the ideal-power torsion functor and a specified projective-resolution completion; the infinitely generated supporting ideals in the examples are not being treated as finitely generated weakly proregular ideals. Localization at the single regular polynomial variable $t$ is used only after the idempotent-ideal Hom functor has been derived.
\end{enumerate}
Thus the paper combines an explicit paired polynomial computation with a concrete separation between degree-zero module conditions and ambient derived behavior. Section~\ref{sec:coordinates} retains the supporting coefficient calculations and examples. Sections~\ref{sec:polynomial} and~\ref{sec:obstruction} establish the comparisons and the main application; Section~\ref{sec:completionexample} treats nonzero completion. The final section records the remaining questions.

\subsection*{Conventions}
All rings are commutative with identity, and all modules are unital. We write $R\Mod$ for all modules, and $R\modf$ for finitely generated modules when $R$ is Noetherian. A finite module means a finitely generated module, not necessarily a finitely presented one. All displayed functor isomorphisms are natural in the coefficient module unless stated otherwise. Negative Ext and Tor groups are zero.

For every ideal $I$ we use
\[
 \Gamma_I(M)=\bigcup_{n\ge1}(0:_M I^n),\qquad
 H_I^q(M)=\RR^q\Gamma_I(M),\qquad
 \Lambda_I(M)=\ilim_n M/I^nM.
\]
Thus $\Gamma_I$ is the small, or ideal-power, torsion functor in the terminology of Rohrer \cite{Rohrer}. For infinitely generated ideals it must be distinguished from large torsion, defined elementwise by $I\subseteq\sqrt{\ann_R(m)}$. No identification with geometric support cohomology is assumed.

For a projective resolution $U_\bullet\to M$, put
\[
 \LL_q\Lambda_I(M)=H_q\bigl(\Lambda_I(U_\bullet)\bigr).
\]
This is the projective-resolution derived functor of the additive completion functor. Independence of resolution follows by chain homotopy comparison; the split horseshoe argument supplies connecting maps, as in \cite[Sections 2.2 and 2.4]{Weibel}. In general $\LL_0\Lambda_I$ need not equal $\Lambda_I$. When $I^2=I$, completion is tensor with $R/I$ and is right exact. Comparisons with other derived completion constructions require additional hypotheses; see \cite{PSY,PSYerr} and \cite[Tag 0BKH]{Stacks}. We do not impose such an identification for infinitely generated ideals.

\section{Reduced modules and relative comparison}
\begin{definition}
An $R$-module $M$ is \emph{$I$-reduced} if $I^2m=0$ implies $Im=0$ for every $m\in M$. It is \emph{$I$-coreduced} if $I^2M=IM$.
\end{definition}
These are the ideal-power conventions of \cite[Definition 2.1]{SsevII}. They need not agree with a definition requiring the corresponding condition for multiplication by every individual element of $I$. For example, take $A=k[\varepsilon]/(\varepsilon^2)$, $R=A\times k$, and $I=A\times0$. Then $I^2=I$, but multiplication by $(\varepsilon,0)$ on $R$ does not satisfy the elementwise kernel condition. We use only the ideal-power definitions above.

\begin{lemma}\label{lem:degreezero}
If $M$ is $I$-reduced, then $(0:_M I^n)=(0:_M I)$ for every $n\ge1$, and
$\Gamma_I(M)\cong\Hom_R(R/I,M)$. If $M$ is $I$-coreduced, then $I^nM=IM$ for every $n\ge1$, and $\Lambda_I(M)\cong M/IM$. These identifications respect the transition maps.
\end{lemma}
\begin{proof}
If $I^nm=0$ for $n\ge2$, apply reducedness to every element of $I^{n-2}m$ to obtain $I^{n-1}m=0$, and then induct. The coreduced assertion follows by multiplying $I^2M=IM$ by powers of $I$. Both defining systems stabilize with their natural maps. Compare \cite[Propositions 2.2--2.3]{SsevI} and \cite[Proposition 2.3]{SsevII}, with the convention just specified.
\end{proof}

\begin{proposition}[Recalled global formulas]\label{prop:global}
For an ideal $I$, the following conditions are equivalent: every module is $I$-reduced; $I^2=I$; every module is $I$-coreduced. If they hold, then on all of $R\Mod$,
\begin{equation}\label{eq:global}
 H_I^q(M)\cong\Ext_R^q(R/I,M),\qquad
 \LL_q\Lambda_I(M)\cong\Tor_q^R(R/I,M).
\end{equation}
For $I\ne R$, the suprema of nonvanishing degrees over all coefficient modules satisfy
\begin{equation}\label{eq:dimensions}
 \cdim_R(I)=\pd_R(R/I),\qquad
 \hdim_R(I)=\fd_R(R/I).
\end{equation}
\end{proposition}
\begin{proof}
The equivalence is \cite[Theorem 3.3]{SsevII}. Directly, the tests $R/I^2$ and $R$ prove the respective implications to idempotence. If $I^2=I$, the functors on the entire module category are $\Hom_R(R/I,-)$ and $(R/I)\otimes_R-$, so deriving gives \eqref{eq:global}. Formula \eqref{eq:dimensions} is the standard characterization of projective and flat dimension by Ext and Tor \cite[Chapter 4]{Weibel}.
\end{proof}
For every ideal, exactness of filtered colimits applied to an injective resolution gives
\begin{equation}\label{eq:extlimit}
 H_I^q(M)\cong\colim_n\Ext_R^q(R/I^n,M).
\end{equation}
See \cite[Tag 00DB]{Stacks} for the exactness used here. For an idempotent ideal, positive torsion vanishes for all coefficients precisely when $R/I$ is projective, and positive completion vanishes for all coefficients precisely when $R/I$ is flat. The first condition is equivalent to $I=eR$ for an idempotent $e$. A finitely generated idempotent ideal always has this form \cite[Tags 00EH, 05KK]{Stacks}.

Let $j:\mathcal E\hookrightarrow R\Mod$ be a full exact inclusion of an abelian category with enough injectives. Define
\[
 H_{I,\mathcal E}^q(M)=\RR^q_{\mathcal E}(\Gamma_I\circ j)(M).
\]
The target is $R\Mod$; no closure under torsion submodules is required.
\begin{proposition}[Standard acyclic-resolution criterion]\label{prop:comparison}
There is a natural map $H_{I,\mathcal E}^q(M)\to H_I^q(jM)$ extending the degree-zero identity. It is an isomorphism for all $q,M$ if and only if every injective object $E$ of $\mathcal E$ satisfies $H_I^q(jE)=0$ for $q>0$.
\end{proposition}
\begin{proof}
Universality of the relative derived cohomological delta functor supplies the map. If internal injectives are ambient acyclic, an internal injective resolution is an acyclic resolution for the ambient functor. Conversely, evaluate an assumed comparison isomorphism on internal injectives. These are the standard arguments of \cite[Sections 2.4--2.5]{Weibel}; see also \cite[Proposition 2.1]{MS}.
\end{proof}
For an ambient additive right exact functor $F$ and enough internal projectives, the homological comparison has direction
\[
 (\LL_q^RF)(jM)\longrightarrow\LL_q^{\mathcal E}(F\circ j)(M),
\]
and is invertible if internal projectives are ambient $F$-acyclic \cite[Section 2.4]{Weibel}. We do not apply this right exact criterion to arbitrary ambient completion.

\begin{proposition}\label{prop:quotient}
Let $K\subseteq R$, $S=R/K$, $J=(I+K)/K$, and $\mathcal E_K=\{M:KM=0\}\simeq S\Mod$. The inclusion is exact, and $\mathcal E_K$ has enough injectives and projectives. All its objects are $I$-reduced, equivalently all are $I$-coreduced, if and only if $J^2=J$. In that case
\[
 H_{I,\mathcal E_K}^q(M)\cong\Ext_S^q(S/J,M),\qquad
 \LL_q^{\mathcal E_K}(\Lambda_I\circ j)(M)\cong\Tor_q^S(S/J,M).
\]
\end{proposition}
\begin{proof}
The powers of $I$ act through the powers of $J$. Apply Proposition~\ref{prop:global} over $S$. Derivation takes place on the entire category $S\Mod$, where the restricted functors are Hom and tensor. The torsion statement also appears in \cite[Proposition 2.2]{MS}.
\end{proof}
\begin{remark}
The category $\mathcal E_K$ is Serre if and only if $K^2=K$. Indeed, an extension of modules killed by $K$ is killed by $K^2$, proving sufficiency. Conversely, extension closure applied to $0\to K/K^2\to R/K^2\to R/K\to0$ forces $K/K^2=0$. An exact quotient subcategory need not therefore be Serre: modules killed by $t$ in $S[t]\Mod$ give an example when $S\ne0$. This distinction matters when applying the Serre-category results of \cite{SsevII}.
\end{remark}

\section{Coordinate calculations and coefficientwise vanishing}\label{sec:coordinates}
Throughout this section $I\subsetneq R$ is generated by a finite or countable family of idempotents. A finite family can be padded with zeros. In particular $I^2=I$, so every $R$-module satisfies
\[
 (0:_M I^2)=(0:_M I),\qquad I^2M=IM.
\]
Thus every coefficient is $I$-reduced and $I$-coreduced. The purpose of the section is to record classical computational tools and their elementary consequences, and to exhibit how positive derived torsion can behave despite these degree-zero conditions. Reducedness and coreducedness impose no additional coefficient restriction in this setting. The same formulas will be used in the later polynomial examples, where the supporting ideal is non-idempotent.

\begin{proposition}[Coordinate formula]\label{prop:coordinate}
Choose idempotent generators $e_1,e_2,\ldots$, and put
\[
 u_0=0,\qquad u_n=1-\prod_{i=1}^n(1-e_i),\qquad d_n=u_n-u_{n-1}.
\]
Then the $d_n$ are pairwise orthogonal idempotents and $I=\bigoplus_{n\ge1}d_nR$. The ideal $I$ is projective, $R/I$ is flat, and for every module $M$ there is a natural exact sequence
\begin{equation}\label{eq:coordinate}
 0\longrightarrow(0:_M I)\longrightarrow M
 \xrightarrow{\Delta_M}\prod_{n\ge1}d_nM
 \longrightarrow H_I^1(M)\longrightarrow0,
 \qquad \Delta_M(m)=(d_nm)_n.
\end{equation}
Moreover,
\[
 H_I^q(M)=0\ (q\ge2),\qquad \LL_q\Lambda_I(M)=0\ (q>0),
\]
and
\[
 \hdim_R(I)=0,\qquad
 \cdim_R(I)=\begin{cases}0,&I\text{ is finitely generated},\\1,&I\text{ is not finitely generated}.\end{cases}
\]
\end{proposition}
\begin{proof}
We have $u_nu_{n-1}=u_{n-1}$, $u_nR=(e_1,\ldots,e_n)$, and $d_n=u_n(1-u_{n-1})$. Thus each $d_n$ is idempotent and is orthogonal to earlier terms. Since $u_n=\sum_{i\le n}d_i$ and $I=\bigcup_nu_nR$, the direct-sum description follows. Each $d_nR$ is a summand of $R$, proving projectivity; this is the classical countable-projectivity result \cite[Tag 05WH]{Stacks} with its explicit generators.

Also $R/I\cong\colim_nR/u_nR$ is a filtered colimit of projective modules, hence is flat \cite[Tags 00DB, 058G]{Stacks}. Applying $\Hom_R(-,M)$ to $0\to I\to R\to R/I\to0$ identifies $\Hom_R(I,M)$ with $\prod_nd_nM$ and gives \eqref{eq:coordinate}, using \eqref{eq:global}. The same length-one projective resolution and flatness give the vanishings. Finally $\pd_R(R/I)\le1$ is zero precisely when $I=eR$ for an idempotent $e$. For the given direct sum this is equivalent to finite generation, or equivalently to only finitely many nonzero summands. Apply \eqref{eq:dimensions}.
\end{proof}

\begin{theorem}[Finite-coefficient vanishing]\label{thm:finite}
Under the hypotheses of Proposition~\ref{prop:coordinate}, the following conditions are equivalent:
\begin{enumerate}[label=(\roman*)]
\item $\rho:R\to\prod_nd_nR$, $r\mapsto(d_nr)_n$, is surjective;
\item $H_I^1(R)=0$;
\item $H_I^q(M)=0$ for every finitely generated $R$-module $M$ and every $q>0$.
\end{enumerate}
If these conditions hold and $I$ is not finitely generated, then
\begin{equation}\label{eq:cyclic}
 H_I^1(I)\cong R/(I+\ann_R I)\ne0.
\end{equation}
Thus vanishing on every finite coefficient need not imply vanishing on the projective coefficient $I$.
\end{theorem}
\begin{proof}
The equivalence of (i) and (ii) is \eqref{eq:coordinate} with $M=R$, and (iii) implies (ii). Suppose (i), and let $m_1,\ldots,m_r$ generate $M$. Given $z_n\in d_nM$, choose $a_{n\ell}\in d_nR$ with $z_n=\sum_{\ell=1}^ra_{n\ell}m_\ell$. These coefficients are obtained by multiplying any expression for $z_n$ by $d_n$. For each $\ell$, surjectivity of $\rho$ gives $a_\ell\in R$ with $d_na_\ell=a_{n\ell}$ for every $n$. Then $m=\sum_\ell a_\ell m_\ell$ satisfies $d_nm=z_n$. Hence $\Delta_M$ is surjective. Higher vanishings follow from Proposition~\ref{prop:coordinate}.

For the coefficient $I$, the image of $\Delta_I$ is exactly $\bigoplus_nd_nR$, so
\begin{equation}\label{eq:productsum}
 H_I^1(I)\cong\frac{\prod_nd_nR}{\bigoplus_nd_nR}.
\end{equation}
The tuple $(d_n)_n$ is not in the direct sum if infinitely many $d_n$ are nonzero; hence the quotient is nonzero when $I$ is not finitely generated. Since $\ker\rho=\ann_R I$ and $\rho(I)=\bigoplus_nd_nR$, the inverse image of the direct sum is $I+\ann_R I$. This proves \eqref{eq:cyclic}.
\end{proof}

\paragraph{The standard Ext argument behind the finite-coefficient test.}
For any $R$-module $A$ with $\pd_R A\le1$, one has
\[
 \Ext_R^1(A,R)=0
 \quad\Longleftrightarrow\quad
 \Ext_R^1(A,M)=0\text{ for every finitely generated }M.
\]
Indeed, choose $0\to K\to R^r\to M\to0$. The long exact Ext sequence contains
\[
 \Ext_R^1(A,R^r)\longrightarrow\Ext_R^1(A,M)
 \longrightarrow\Ext_R^2(A,K),
\]
and the outer terms vanish under the stated hypothesis. The converse follows by taking $M=R$. This uses only the usual long exact sequence and projective-dimension criterion \cite[Chapters 2 and 4]{Weibel}. Taking $A=R/I$ explains the finite-coefficient equivalence in Theorem~\ref{thm:finite}; its coordinate formulation follows from \eqref{eq:coordinate}.

The following sum criterion is likewise a consequence of the explicit projective presentation. For the general theory of preservation of sums by defect and Ext functors, including criteria expressed through projective presentations, see \cite{BZ}. We give the short argument specialized to the orthogonal decomposition, which also identifies a concrete family witnessing failure.

\begin{theorem}[Countable direct sums]\label{thm:sums}
Under the hypotheses of Proposition~\ref{prop:coordinate}, the following are equivalent:
\begin{enumerate}[label=(\roman*)]
\item $I$ is finitely generated;
\item $H_I^1$ is the zero functor;
\item $H_I^1$ commutes with countable direct sums;
\item the canonical map $\bigoplus_nH_I^1(d_nR)\to H_I^1(\bigoplus_nd_nR)$ is an isomorphism.
\end{enumerate}
\end{theorem}
\begin{proof}
If $I$ is finitely generated, it is generated by an idempotent, so the functor is zero. Thus (i) implies (ii), which implies (iii), and (iii) implies (iv). For each $n$, the only possibly nonzero coordinate of $\Delta_{d_nR}$ is the $n$th, and this coordinate is the identity; hence $H_I^1(d_nR)=0$. If (iv) holds, then $H_I^1(I)=0$, and \eqref{eq:productsum} forces only finitely many nonzero $d_n$. Thus (iv) implies (i).
\end{proof}
Theorem~\ref{thm:sums} does not require surjectivity of $\rho$. In particular, preservation of countable direct sums is a stronger demand than vanishing on all finite coefficients.

\subsection{Two contrasting coefficient calculations}
\begin{example}[Product of fields]\label{ex:product}
Let $R=\prod_{n\ge1}k$, $I=\bigoplus_{n\ge1}k$, and $T=R/I$, for a field $k$. The coordinate idempotents give an isomorphism $\rho:R\to\prod_ne_nR$ and $\ann_R I=0$. Consequently every finite coefficient has zero positive torsion, whereas
\[
 H_I^q(I)\cong\begin{cases}T,&q=1,\\0,&q\ne1.\end{cases}
\]
Indeed $(0:_I I)=0$ and \eqref{eq:productsum} gives degree one. Also $\cdim_R(I)=1$ and $\hdim_R(I)=0$. This is the cyclic calculation of \cite[Example 5.4]{MS}, used below as an input to the polynomial construction.
\end{example}

\begin{example}[A rank-one free coefficient with nonzero torsion]\label{ex:eventual}
Let
\[
 R_{\mathrm{ec}}=\{(a_n)\in k^{\NN}:(a_n)\text{ is eventually constant}\},
 \qquad I=k^{(\NN)}.
\]
This is a unital subring of $k^{\NN}$ and $I$ is its ideal of finitely supported sequences. Since $I=\bigoplus_ne_nR_{\mathrm{ec}}$ and $e_nR_{\mathrm{ec}}\cong k$, the coordinate map is the inclusion $R_{\mathrm{ec}}\hookrightarrow k^{\NN}$. Hence
\[
 H_I^1(R_{\mathrm{ec}})\cong k^{\NN}/R_{\mathrm{ec}}\ne0,
 \qquad H_I^1(I)\cong k^{\NN}/k^{(\NN)}.
\]
The alternating sequence $0,1,0,1,\ldots$ is not eventually constant, proving the first nonvanishing over every field. Thus a free coefficient of rank one can have nonzero first derived torsion when coordinate surjectivity fails. As before, all torsion in degrees at least two and all positive completion vanish. This illustrates finite-coefficient nonvanishing outside the coordinate-surjectivity hypothesis; the calculation is a direct application of Proposition~\ref{prop:coordinate}.
\end{example}

\section{Polynomial comparison and bounded nilpotence}\label{sec:polynomial}
Let $S$ be a commutative ring, $J\subseteq S$ an idempotent ideal, and $T=S/J$. In this section set
\begin{equation}\label{eq:polynomialsetup}
\begin{gathered}
 R=S[t],\qquad I=JR+(t),\qquad A_N=S[t]/(t^N),\\
 \mathcal E_N=\{M\in R\Mod:t^NM=0\}\simeq A_N\Mod.
\end{gathered}
\end{equation}
The integer $N\ge1$ is fixed when defining $\mathcal E_N$.

\begin{lemma}[Compatible cone resolutions]\label{lem:cones}
For every $n\ge1$,
\[
 I^n=JR+(t^n),\qquad R/I^n\cong T[t]/(t^n).
\]
If $P_\bullet\to T$ is an $S$-projective resolution and $Q_\bullet=R\otimes_SP_\bullet$, then a projective resolution of $R/I^n$ is
\[
 C_i^{(n)}=Q_i\oplus Q_{i-1},\qquad
 d^{(n)}(a,b)=(d_Qa+t^nb,-d_Qb),\qquad Q_{-1}=0.
\]
The canonical quotient $R/I^{n+1}\to R/I^n$ is lifted by
\begin{equation}\label{eq:conemap}
 f_n:C_\bullet^{(n+1)}\longrightarrow C_\bullet^{(n)},
 \qquad f_n(a,b)=(a,tb).
\end{equation}
\end{lemma}
\begin{proof}
Expansion of powers and $J^n=J$ give the formula for $I^n$. The ring $R$ is free over $S$, so $Q_\bullet$ is a projective resolution of $T[t]$. Multiplication by $t^n$ is injective on $T[t]$. The cone homology sequence therefore proves that the displayed cone resolves its cokernel $T[t]/(t^n)$ \cite[Section 1.5]{Weibel}. The identity $d^{(n)}f_n=f_nd^{(n+1)}$ follows directly, and the augmentation induces the quotient map.
\end{proof}

\begin{theorem}[Bounded-nilpotence comparison]\label{thm:bounded}
In \eqref{eq:polynomialsetup}, let $t^NM=0$. For every $q\ge0$ there are natural isomorphisms
\begin{align}
 H_I^q(M)&\cong\Ext_S^q(T,M)\cong H_{I,\mathcal E_N}^q(M),\label{eq:boundedtorsion}\\
 \LL_q\Lambda_I(M)&\cong\Tor_q^S(T,M)
 \cong\LL_q^{\mathcal E_N}(\Lambda_I\circ j)(M).
 \label{eq:boundedcompletion}
\end{align}
More precisely, for $n\ge N$,
\begin{align}
 \Ext_R^q(R/I^n,M)&\cong\Ext_S^q(T,M)\oplus\Ext_S^{q-1}(T,M),\label{eq:extstages}\\
 \Tor_q^R(R/I^n,M)&\cong\Tor_q^S(T,M)\oplus\Tor_{q-1}^S(T,M).
 \label{eq:torstages}
\end{align}
The Ext transitions run from stage $n$ to stage $n+1$, and the Tor transitions run from stage $n+1$ to stage $n$; on both displayed decompositions their formula is $(a,b)\mapsto(a,tb)$. Their $N$-fold composites kill the shifted summands. If $T\ne0$, then $I^2\ne I$.
\end{theorem}
\begin{proof}
For $n\ge N$, the $t^n$ maps vanish after applying either $\Hom_R(-,M)$ or $-\otimes_RM$ to the cone of Lemma~\ref{lem:cones}. By adjunction the resulting two complexes are the direct sums of $\Hom_S(P_\bullet,M)$ and its cohomological degree-one shift, or $P_\bullet\otimes_SM$ and its homological degree-one shift, respectively. This proves \eqref{eq:extstages} and \eqref{eq:torstages}. The maps \eqref{eq:conemap} induce $(a,b)\mapsto(a,tb)$, where $t$ acts through its endomorphism of $M$. Balancing of Tor identifies this calculation with the homology obtained from a projective resolution of $M$ \cite[Section 2.7]{Weibel}.

The shifted direct system has zero $N$-fold maps, so its colimit is zero. Formula \eqref{eq:extlimit} now gives the first isomorphism of \eqref{eq:boundedtorsion}.

For completion choose an ambient projective resolution $U_\bullet\to M$. The tower $U_\bullet/I^nU_\bullet$ is degreewise surjective. The Milnor exact sequence \cite[Section 3.5]{Weibel} gives
\begin{equation}\label{eq:milnor}
 0\longrightarrow\ilim_n{}^{1}\Tor_{q+1}^R(R/I^n,M)
 \longrightarrow H_q\bigl(\ilim_nU_\bullet/I^nU_\bullet\bigr)
 \longrightarrow\ilim_n\Tor_q^R(R/I^n,M)\longrightarrow0.
\end{equation}
The homology tower is, from stage $N$ onward, the constant first summand and a second summand with $N$-fold zero maps. For every fixed stage, all sufficiently distant images are exactly the first summand, so the tower is Mittag--Leffler. Its inverse limit is $\Tor_q^S(T,M)$ and its first derived inverse limit is zero. This proves the first isomorphism of \eqref{eq:boundedcompletion}. Notice that degreewise surjectivity of the complex tower and the Mittag--Leffler property of its homology tower are distinct steps.

On $\mathcal E_N$, the identities $I^nM=JM$ for $n\ge N$ and
\[
 \Gamma_I(M)=(0:_M J),\qquad \Lambda_I(M)=M/JM
\]
hold. For torsion, if $I^nm=0$, then $Jm=0$ since $J\subseteq I^n$; conversely $Jm=0$ and $t^NM=0$ imply $I^Nm=0$. Thus the restricted functors are $\Hom_S(T,-)$ and $T\otimes_S-$. Since $A_N$ is a finite free $S$-module, restriction of scalars sends $A_N$-projectives to $S$-projectives. It also sends $A_N$-injectives to $S$-injectives: its left adjoint $A_N\otimes_S-$ is exact. Internal resolutions consequently compute the same Ext and Tor over $S$. This proves the remaining isomorphisms. In the torsion case, an internal injective is $S$-injective, so the already computed ambient groups vanish on it in positive degrees. Proposition~\ref{prop:comparison} therefore identifies the canonical map $H_{I,\mathcal E_N}^q(M)\to H_I^q(M)$ as an isomorphism. The completion identities above are obtained from the explicit projective resolutions, without invoking a right exact comparison theorem for the unrestricted ambient completion functor.

Finally, reduction modulo $JR$ sends $I$ to $(t)$ in $T[t]$, whose square is different from itself if $T\ne0$.
\end{proof}

\begin{remark}\label{rem:reducedboundary}
When $N=1$, all objects of $\mathcal E_1$ are $I$-reduced and $I$-coreduced by Proposition~\ref{prop:quotient}. This is the reduced-module setting of the original polynomial comparison. Theorem~\ref{thm:bounded} extends both calculations to bounded nilpotence. For $N>1$, $\mathcal E_N$ is not in general a category of reduced or coreduced modules: take $S=k$, $J=0$, and $M=k[t]/(t^N)$. For bounded coefficients the comparison is therefore more general than the reduced-module setting.
\end{remark}

\begin{corollary}\label{cor:dimensions}
For every fixed $N\ge1$ and $T\ne0$, the suprema over $\mathcal E_N$ of the nonvanishing degrees of ambient torsion and projective-derived completion are $\pd_S(T)$ and $\fd_S(T)$, respectively. They are also the corresponding relative dimensions on $\mathcal E_N$.
\end{corollary}
\begin{proof}
Theorem~\ref{thm:bounded} gives the upper bounds. Every $S$-module belongs to $\mathcal E_N$ after letting $t$ act by zero, so the usual Ext and Tor tests give equality \cite[Chapter 4]{Weibel}.
\end{proof}

\begin{corollary}\label{cor:liftfinite}
Suppose $J\subsetneq S$ is generated by countably many idempotents with orthogonal form $J=\bigoplus_nd_nS$, and $S\to\prod_nd_nS$ is surjective. For each fixed $N$, every finitely generated $A_N$-module has zero positive ambient and relative $I$-torsion, and every object of $\mathcal E_N$ has zero positive ambient and relative projective-derived completion. If $J$ is not finitely generated, then the coefficient $J$ with $tJ=0$ has
\[
 H_I^1(J)\cong S/(J+\ann_SJ)\ne0.
\]
Moreover, on $\mathcal E_N$ the functor $H_I^1$ commutes with countable direct sums if and only if $J$ is finitely generated.
\end{corollary}
\begin{proof}
A finite $A_N$-module is finite over $S$, since $A_N$ is finite free over $S$. Apply Theorem~\ref{thm:finite} over $S$ and then Theorem~\ref{thm:bounded}. Flatness of $S/J$ proves completion vanishing. The direct-sum assertion follows from Theorem~\ref{thm:sums}, testing necessity on the modules $d_nS$ with zero $t$-action and using $\Ext_S^1(S/J,-)=0$ for sufficiency when $J$ is finitely generated.
\end{proof}

\section{The obstruction beyond bounded nilpotence}\label{sec:obstruction}
We retain $R=S[t]$, $I=JR+(t)$, $J^2=J$, and $T=S/J$, but allow an arbitrary $R$-module $M$. Multiplication by $t$ on $M$ makes each $\Ext_S^q(T,M)$ an $S[t]$-module. For any such module $E$, write $E[t^{-1}]$ for its localization at the powers of $t$, and $\lambda_E:E\to E[t^{-1}]$ for the localization map.

\begin{theorem}[Polynomial torsion and its comparison map]\label{thm:fiber}
Let $P_\bullet\to T$ be an $S$-projective resolution and set $A^\bullet=\Hom_S(P_\bullet,M)$, with its induced $t$-action. The complex with terms and differential
\[
 K^q=A^q\oplus A^{q-1}[t^{-1}],\qquad
 d_K(a,b)=(d_Aa,\lambda_A(a)-d_Ab)
\]
computes $H_I^q(M)$. Equivalently, it is a model for the mapping fiber of $A^\bullet\to A^\bullet[t^{-1}]$. If $E^q=\Ext_S^q(T,M)$, there are natural short exact sequences
\begin{equation}\label{eq:obstruction}
 0\longrightarrow\coker(\lambda_{E^{q-1}})
 \longrightarrow H_I^q(M)
 \xrightarrow{\alpha_q}\ker(\lambda_{E^q})
 \longrightarrow0.
\end{equation}
Here $E^{-1}=0$. The induced map $H_I^q(M)\to E^q$ is an isomorphism if and only if $\lambda_{E^{q-1}}$ is surjective and $\lambda_{E^q}$ is the zero map. These maps are isomorphisms for every $q$ if and only if every $E^q$ is $t$-power torsion.
\end{theorem}
\begin{proof}
The cone resolutions of Lemma~\ref{lem:cones} identify $\Hom_R(C_\bullet^{(n)},M)$ with the cochain complex
\[
 D_n^q=A^q\oplus A^{q-1},\qquad
 d_n(a,b)=(d_Aa,t^na-d_Ab).
\]
The direct-system map is $(a,b)\mapsto(a,tb)$. The first summands therefore have colimit $A^\bullet$, and the second summands have colimit $A^\bullet[t^{-1}]$. Explicitly, a second component $b$ at stage $n$ is sent to $b/t^n$. Under these identifications the cross map $t^n$ becomes localization, so $\colim_nD_n^\bullet=K^\bullet$. Exactness of filtered colimits and \eqref{eq:extlimit} identify its cohomology with $H_I^q(M)$.

The projection to $A^\bullet$ and inclusion of the shifted localized complex yield the mapping-fiber long exact sequence. Localization is exact, so $H^q(A^\bullet[t^{-1}])=E^q[t^{-1}]$. The corresponding kernel--cokernel pieces are \eqref{eq:obstruction}. These are applications of the usual cone sequence and exactness of filtered colimits \cite[Section 1.5]{Weibel}, \cite[Tag 00DB]{Stacks}.

The map to $E^q$ has kernel $\coker\lambda_{E^{q-1}}$ and image $\ker\lambda_{E^q}$, which proves the degreewise criterion. If $E^q$ is $t$-power torsion, then $E^q[t^{-1}]=0$, and all maps are isomorphisms. Conversely, isomorphisms in all degrees force $\lambda_{E^q}$ to be zero for every $q$. A zero localization map has zero localized target, since every localized element is a fraction of an element in its image. Hence every $E^q$ is $t$-power torsion.
\end{proof}

\paragraph{Derived-composition interpretation and the existing comparison sequence.}
Put $F=\Hom_R(R/JR,-)=\Hom_R(T[t],-)$ and $G=\Gamma_{(t)}$. Since $I^n=JR+(t^n)$, one has $\Gamma_I=G\circ F$. For an injective $R$-module $E$, the module $F(E)$ is injective over $T[t]$: it is the right adjoint image for the exact restriction-of-scalars functor. Multiplication by $t$ on $T[t]$ is injective, so applying $\Hom_{T[t]}(-,F(E))$ shows that $t:F(E)\to F(E)$ is surjective.

For completeness, the two-term complex $N\to N[t^{-1}]$ computes $\RR\Gamma_{(t)}(N)$ over $R$. Indeed, $t$ is a nonzerodivisor on $R$, so multiplication by $t$ is surjective on every injective $R$-module. Its localization map is then surjective as well. Exactness of localization makes the cohomology of this two-term complex a cohomological delta functor with degree zero $\Gamma_{(t)}$ and zero positive values on injectives, which proves the claim. The same surjectivity argument makes $F(E)$ acyclic for $G$.

The derived-composition theorem and Grothendieck spectral sequence \cite[Tags 015M, 015N]{Stacks} consequently give
\[
 \RR\Gamma_I\cong\RR\Gamma_{(t)}\circ\RR F,\qquad
 E_2^{p,q}=H_{(t)}^p\bigl(\Ext_R^q(T[t],M)\bigr)
 \Longrightarrow H_I^{p+q}(M).
\]
Base change of the projective resolution $P_\bullet\to T$ identifies $\Ext_R^q(T[t],M)$ with $\Ext_S^q(T,M)$. Only the columns $p=0,1$ occur. Their filtration yields exactly \eqref{eq:obstruction}, while the cone calculation above records the map to $E^q$ explicitly.

This identifies the formula with the established comparison-sequence framework. In particular, \cite[Proposition 3.4]{QR} obtains
\[
 0\to H_{(b)}^1 H_{\mathfrak a}^{q-1}(M)
 \to H_{\mathfrak a+(b)}^q(M)
 \to\Gamma_{(b)}H_{\mathfrak a}^q(M)\to0
\]
under injective-preservation assumptions. Here the required acyclicity for the particular pair $F,G$ has been verified directly, without assuming preservation of ambient injectives by both torsion functors. The use of a two-term localization complex concerns $(t)$; it does not assert a \v{C}ech model for the infinitely generated ideal $I$.

\begin{remark}
Theorem~\ref{thm:fiber} does not identify $\Ext_S^q(T,M)[t^{-1}]$ with $\Ext_S^q(T,M[t^{-1}])$. Such an interchange would require extra hypotheses. The localization in the formula is applied after forming the Hom complex. This is precisely where unbounded locally nilpotent coefficients can behave differently from bounded ones.
\end{remark}

\begin{corollary}\label{cor:localcriterion}
If multiplication by $t$ is locally nilpotent on every $\Ext_S^q(T,M)$, then $H_I^q(M)\cong\Ext_S^q(T,M)$ for all $q$. In particular this applies when $t^NM=0$ for some fixed $N$.
\end{corollary}
\begin{proof}
Apply Theorem~\ref{thm:fiber}. If $t^NM=0$, the same power kills every term of $\Hom_S(P_\bullet,M)$ and hence every Ext group.
\end{proof}
Local nilpotence on $M$ alone is insufficient. The following construction proves the coefficient assertion in Theorem~\ref{thm:main}.

\begin{theorem}[A reduced and coreduced coefficient with second torsion]\label{ex:localnil}
Let $S=\prod_{n\ge1}k$, $J=\bigoplus_{n\ge1}k$, $T=S/J$, $R=S[t]$, and $I=JR+(t)$. Let
\[
 V_n=k[t]/(t^n),\qquad M=\bigoplus_{n\ge1}V_n,
\]
where $S$ acts on $V_n$ through its $n$th coordinate and $t$ acts by the polynomial action. Then $I^2\ne I$, $M$ is $I$-reduced and $I$-coreduced, and
\[
 \Gamma_I(M)=\Lambda_I(M)=0,\qquad H_I^2(M)\ne0.
\]
Writing $B=(\prod_nV_n)/(\bigoplus_nV_n)$, one has
\[
 H_I^1(M)\cong\ker(B\to B[t^{-1}]),\qquad
 H_I^2(M)\cong\coker(B\to B[t^{-1}]),
\]
with $H_I^q(M)=0$ for $q=0$ and $q\ge3$. Each $M_N=\bigoplus_{1\le n\le N}V_n$ is finite, $I$-reduced and $I$-coreduced, and has $H_I^q(M_N)=0$ for $q>0$, although $M=\colim_NM_N$ has nonzero degree two.
\end{theorem}
\begin{proof}
Let $e_n\in S$ denote the $n$th coordinate idempotent. Every element of $M$ is killed by some power of $t$, but no common power kills $M$.

We first verify the module-theoretic conditions. If $m\in M$ has support in a finite set $F\subseteq\NN$, the idempotent $e_F=\sum_{n\in F}e_n\in J$ satisfies $e_Fm=m$. Hence $JM=M$, and therefore
\[
 IM=M=I^2M.
\]
Thus $M$ is $I$-coreduced. If $Jm=0$, then $e_nm=0$ for every $n$, so every coordinate of $m$ is zero. Consequently $(0:_M J)=0$. Since $JR\subseteq I^2\subseteq I$, it follows that
\[
 (0:_M I)=(0:_M I^2)=0,
\]
which proves that $M$ is also $I$-reduced. Lemma~\ref{lem:degreezero} now gives
\[
 \Gamma_I(M)=0,\qquad \Lambda_I(M)=0.
\]
The supporting ideal is nevertheless non-idempotent, since its image in $T[t]$ is $(t)$ and $T\ne0$.

We next compute the ambient derived torsion. Proposition~\ref{prop:coordinate} over $S$ gives
\[
 E^0=\Hom_S(T,M)=0,\qquad
 E^1=\Ext_S^1(T,M)\cong B:=\frac{\prod_nV_n}{\bigoplus_nV_n},
 \qquad E^q=0\ (q\ge2).
\]
Indeed the coordinate map on $M$ is the direct-sum inclusion into $\prod_nV_n$. Theorem~\ref{thm:fiber} yields
\begin{equation}\label{eq:localnilgroups}
 H_I^1(M)\cong\ker(B\to B[t^{-1}]),\qquad
 H_I^2(M)\cong\coker(B\to B[t^{-1}]),
\end{equation}
and $H_I^q(M)=0$ for $q=0$ and $q\ge3$.

Both the failure of comparison and nonzero degree two can be checked explicitly. Let $b\in B$ be the class of $(1\bmod t^n)_n$. For every $r\ge0$, its multiple $t^rb$ has nonzero coordinates for every $n>r$. Hence $t^rb\ne0$ in $B$, so $b$ is not in the image of the comparison $H_I^1(M)\to B$.

The class $b/t$ is not in the image of $B\to B[t^{-1}]$. Otherwise $b/t=c/1$ for some $c\in B$, so $t^r(b-tc)=0$ in $B$ for some $r\ge0$. Choose a representative $(c_n)_n$ of $c$. For every $n>r$, the polynomial
\[
 t^r(1-tc_n)\in k[t]/(t^n)
\]
has nonzero coefficient $1$ in degree $r$, so is nonzero. This contradicts vanishing in the product modulo the direct sum. Thus
\[
 H_I^2(M)\ne0.
\]
This establishes a simultaneous reduced and coreduced coefficient for a non-idempotent ideal with
\[
 \Gamma_I(M)=\Lambda_I(M)=0
 \qquad\text{but}\qquad H_I^2(M)\ne0.
\]
This is ambient derived torsion; the degree-zero identities are being tested on the coefficient rather than derived inside a category of all reduced and coreduced modules.

The same coordinate-idempotent argument proves reducedness and coreducedness of each finite partial sum. Each $M_N=\bigoplus_{1\le n\le N}V_n$ is a finite $S$-module killed by $t^N$, so $H_I^q(M_N)=0$ for $q>0$ by Example~\ref{ex:product} and Theorem~\ref{thm:bounded}. Since $M=\colim_NM_N$, this also gives an explicit failure of $H_I^2$ to commute with this filtered colimit. There is no contradiction with Proposition~\ref{prop:coordinate}: the ideal here is $JR+(t)$ and is not idempotent, and the coefficient is not killed by a fixed power of $t$.
\end{proof}

\begin{corollary}[Unrestricted and bounded torsion dimensions]\label{cor:dimensionjump}
For $S=\prod_{n\ge1}k$, $J=\bigoplus_{n\ge1}k$, $R=S[t]$, and $I=JR+(t)$, one has
\[
 \cdim_S(J)=1,\qquad \cdim_R(I)=2,
 \qquad \sup_{M\in\mathcal E_N}\{q:H_I^q(M)\ne0\}=1
 \quad(N\ge1).
\]
In the supremum, all coefficients in $\mathcal E_N$ are allowed. If it is restricted to finitely generated $A_N$-modules, every positive group vanishes.
\end{corollary}
\begin{proof}
Proposition~\ref{prop:coordinate} gives $\pd_S(T)=1$, so $E^q=\Ext_S^q(T,M)=0$ for $q\ge2$ and every $R$-module $M$. The obstruction sequence therefore gives $H_I^q(M)=0$ for $q\ge3$. Theorem~\ref{ex:localnil} gives a coefficient with $H_I^2(M)\ne0$, proving the unrestricted value two. The bounded value follows from Corollary~\ref{cor:dimensions}, or is detected by $J$ with zero $t$-action. Finite-coefficient vanishing is Corollary~\ref{cor:liftfinite}.
\end{proof}

\section{A nonzero derived completion on a cyclic coefficient}\label{sec:completionexample}
The countably idempotent-generated class has flat quotient and therefore zero positive derived completion. Idempotence of an ideal alone does not imply this flatness. We now give an explicit calculation outside that class.

\begin{example}[An idempotent maximal ideal in a valuation domain]\label{ex:valuation}
Fix a field $k$ and compatible indeterminates $x_n=x^{1/n!}$, with $x_n=x_{n+1}^{n+1}$. Inside the union of their fraction fields put
\[
 V=\bigcup_{n\ge1} k[x_n]_{(x_n)},\qquad
 \mathfrak m=\bigcup_{n\ge1}x_nV.
\]
The indicated local inclusions are well defined: a denominator with nonzero constant term remains a unit at every later stage. The union is a local domain with maximal ideal $\mathfrak m$ and residue field $k$. It is a valuation domain. Indeed, every nonzero element of its fraction field belongs to the fraction field of some $k[x_n]_{(x_n)}$, which is a discrete valuation ring; either the element or its inverse lies in that ring, hence in $V$.

Every $a\in\mathfrak m$ lies in $x_nV$ for some $n$. Since
\[
 x_n=x_{n+1}^{n+1}=x_{n+1}\,x_{n+1}^{n},
\]
we have $x_n\in\mathfrak m^2$. Thus $\mathfrak m^2=\mathfrak m$. The ideal is nonzero and proper and cannot be generated by idempotent elements, since a domain has only the idempotents $0$ and $1$.

Put $a=x_1=x$ and $C=V/aV$. Because $a$ is a nonzero element of a domain, a free resolution of this cyclic coefficient is
\[
 0\longrightarrow V\xrightarrow{a}V\longrightarrow C\longrightarrow0.
\]
Tensoring with $V/\mathfrak m\cong k$ makes the displayed multiplication map zero. Consequently
\begin{equation}\label{eq:valuationcompletion}
 \LL_q\Lambda_{\mathfrak m}(C)
 \cong\Tor_q^V(k,C)
 \cong\begin{cases}k,&q=0,1,\\0,&q\ge2.\end{cases}
\end{equation}
Here Proposition~\ref{prop:global} applies because $\mathfrak m^2=\mathfrak m$. In particular $k$ is not flat over $V$, as is already detected by $\Tor_1^V(k,C)\ne0$. Every $V$-module is nevertheless $\mathfrak m$-reduced and $\mathfrak m$-coreduced.

For the polynomial extension $R=V[t]$, $I=\mathfrak mR+(t)$, give $C$ the zero $t$-action. Theorem~\ref{thm:bounded} then gives
\[
 \LL_1\Lambda_I(C)\cong k\ne0,
 \qquad \LL_q\Lambda_I(C)=0\quad(q\ge2),
\]
with the same formulas for relative completion in $\mathcal E_1$. This is a non-idempotent supporting ideal with nonzero first derived completion on a cyclic reduced and coreduced coefficient.
\end{example}

\begin{remark}
The explicit free resolution of $C$ computes the specified projective-derived completion. The nonzero group in \eqref{eq:valuationcompletion} shows why the flat-quotient condition in the earlier vanishing results is substantive.
\end{remark}

\section{Scope and further questions}
The standard coefficientwise consequences recorded in Section~\ref{sec:coordinates} give two different tests. Coordinate surjectivity is equivalent to vanishing on all finite coefficients, whereas finite generation of the ideal is equivalent to preservation of countable direct sums by first torsion. The ring of eventually constant sequences supplies finite-coefficient nonvanishing when coordinate surjectivity fails, while a product of fields gives the opposite finite-coefficient behavior. Polynomial lifting preserves these bounded-coefficient calculations.

The bounded-nilpotence theorem and the mapping-fiber formula also explain the limit of that comparison. Local nilpotence of $t$ on the coefficient alone does not ensure local nilpotence on its Ext groups. Theorem~\ref{ex:localnil} detects the resulting obstruction in degree two on a coefficient that is simultaneously reduced and coreduced for a non-idempotent ideal. It therefore supplies a direct module-theoretic application beyond the automatic reducedness and coreducedness of the idempotent setting. On the completion side, Example~\ref{ex:valuation} gives a cyclic coefficient with nonzero degree one. These examples address the finite-coefficient and nontrivial-completion issues without asserting a general classification.

\begin{question}
For idempotent ideals outside the countably idempotent-generated class, which pairs $(\pd_R(R/I),\fd_R(R/I))$ occur? Can explicit coefficients realizing higher nonzero derived completion be described with comparably short resolutions?
\end{question}
\begin{question}
For the polynomial setup, which additional hypotheses on $S$, $J$, and a locally nilpotent coefficient $M$ force $\Ext_S^q(S/J,M)$ to remain $t$-power torsion? Theorem~\ref{thm:fiber} gives an exact criterion, but a structural classification beyond bounded nilpotence remains to be developed.
\end{question}
\begin{question}
Which nonquotient exact full abelian subcategories of $I$-reduced modules have enough injectives satisfying the ambient acyclicity criterion? Exactness, enough injectives, and ambient acyclicity must be checked separately.
\end{question}
For infinitely generated supporting ideals, comparisons with geometric or sequential torsion and completion require an additional theorem. The calculations here concern precisely the ideal-power torsion and projective-resolution completion fixed in the conventions.

\section*{Acknowledgement}
The author gratefully acknowledges support from the Eastern Africa Algebra Research Group. The author used ChatGPT (OpenAI) to assist with the revision of the manuscript.
\section*{Disclosure statement}
The author declares no conflict of interest.


\begin{thebibliography}{99}
\bibitem{BZ} S. Breaz and J. \v{Z}emli\v{c}ka, The defect functor of a homomorphism and direct unions, \emph{Algebr. Represent. Theory} \textbf{19} (2016), no.~1, 181--208. \href{https://doi.org/10.1007/s10468-015-9569-0}{doi:10.1007/s10468-015-9569-0}.
\bibitem{MS} P. Masaba and D. Ssevviiri, \emph{Relative derived torsion in categories of reduced modules: comparison and non-Noetherian examples},
https://doi.org/10.48550/arXiv.2611.13603.
\bibitem{PSY} M. Porta, L. Shaul and A. Yekutieli, On the homology of completion and torsion, \emph{Algebr. Represent. Theory} \textbf{17} (2014), no.~1, 31--67. \href{https://doi.org/10.1007/s10468-012-9385-8}{doi:10.1007/s10468-012-9385-8}.
\bibitem{PSYerr} M. Porta, L. Shaul and A. Yekutieli, Erratum to: On the homology of completion and torsion, \emph{Algebr. Represent. Theory} \textbf{18} (2015), no.~5, 1401--1405. \href{https://doi.org/10.1007/s10468-015-9557-4}{doi:10.1007/s10468-015-9557-4}.
\bibitem{QR} P.~H. Qu\'y and F. Rohrer, Injective modules and torsion functors, \emph{Comm. Algebra} \textbf{45} (2017), no.~1, 285--298. \href{https://doi.org/10.1080/00927872.2016.1206345}{doi:10.1080/00927872.2016.1206345}.
\bibitem{Rohrer} F. Rohrer, Torsion functors, small or large, \emph{Beitr. Algebra Geom.} \textbf{60} (2019), 233--256. \href{https://doi.org/10.1007/s13366-018-0414-6}{doi:10.1007/s13366-018-0414-6}.
\bibitem{SsevI} D. Ssevviiri, Applications of reduced and coreduced modules I, \emph{Int. Electron. J. Algebra} \textbf{35} (2024), 61--81. \href{https://doi.org/10.24330/ieja.1299587}{doi:10.24330/ieja.1299587}.
\bibitem{SsevII} D. Ssevviiri, \emph{Applications of reduced and coreduced modules II: Radicality of the functor $\Hom_R(R/I,-)$}, accepted for publication as a chapter in the book: Rings and Polynomials: Algebraic, Arithmetic, and Topological Topics in Ring Theory; preprint \href{https://arxiv.org/abs/2306.12871}{arXiv:2306.12871}.
\bibitem{Stacks} The Stacks Project Authors, \emph{The Stacks Project}, \url{https://stacks.math.columbia.edu}. Tags 00DB, 00EH, 015M, 015N, 05KK, 05WH, 058G, and 0BKH.
\bibitem{VY} R. Vyas and A. Yekutieli, Weak proregularity, weak stability, and the noncommutative MGM equivalence, \emph{J. Algebra} \textbf{513} (2018), 265--325. \href{https://doi.org/10.1016/j.jalgebra.2018.07.023}{doi:10.1016/j.jalgebra.2018.07.023}.
\bibitem{Weibel} C.~A. Weibel, \emph{An Introduction to Homological Algebra}, Cambridge Studies in Advanced Mathematics, vol.~38, Cambridge University Press, Cambridge, 1994.
\bibitem{Yek} A. Yekutieli, Weak proregularity, derived completion, adic flatness, and prisms, \emph{J. Algebra} \textbf{583} (2021), 126--152. \href{https://doi.org/10.1016/j.jalgebra.2021.04.033}{doi:10.1016/j.jalgebra.2021.04.033}.
\end{thebibliography}
\end{document}